\documentclass[a4paper,12pt]{article}
\usepackage{graphicx}
\usepackage{amsthm,amsmath, amsfonts}
\usepackage{amsthm,amscd,amsfonts,latexsym,color,epsfig}
\usepackage{amssymb, epic, graphicx}

\title{Boundary-value problem for fractional heat equation involving Caputo-Fabrizio derivative}
\author{Nasser Al-Salti, Erkinjon Karimov and Sebti Kerbal}

\begin{document}
\sloppy

\maketitle

\begin{abstract}
In this work, we consider a number of boundary-value problems for time-fractional heat equation with the recently introduced Caputo-Fabrizio derivative. Using the method of separation of variables, we prove a unique solvability of the stated problems. Moreover, we have found an explicit solution to certain initial value problem for Caputo-Fabrizio fractional order differential equation  by reducing the problem to a Volterra integral equation. Different forms of solution were presented depending on the values of the parameter appeared in the problem.
\end{abstract}

\section{Introduction and Preliminary}

\subsection{Definitions and related works}

Recently, Caputo and Fabrizio introduced a new fractional derivative [1]
\begin{equation}\label{eq1}
_{CF}D_{at}^\alpha f(t)=\frac{1}{1-\alpha}\int\limits_a^t f'(s)e^{-\frac{\alpha}{1-\alpha}(t-s)}ds,
\end{equation}
where the order of the derivative $\alpha\in [0,1]$. In their next work [2], they defined a domain, on which operator (1) is well defined by the set 
$$
W^{\alpha,1}(a,\infty)=\left\{f(t)\in L^1(a,\infty);\, \left(f(t)-f_a(s)\right)e^{-\frac{\alpha}{1-\alpha}(t-s)}\in L^1(a,t)\times L^1(a,\infty)\right\},
$$ 
whose norm is given for $\alpha\neq 1$ by
$$
||f(t)||_{W^{\alpha,1}}=\int\limits_a^\infty|f(t)|dt+\frac{\alpha}{1-\alpha}\int\limits_a^\infty\int\limits_{-\infty}^t |f_a(s)|e^{-\frac{\alpha}{1-\alpha}(t-s)}dsdt,
$$
where $f_a$ denotes the extension of the function $f(t)$ and given by
$$
f_a(t)=f(t),\,\,\,t\geq a,\,\,\,\,\,\,\,\,\,f_a(t)=0,\,\,-\infty<t<a.
$$

The interest of this new fractional derivative is justified by Caputo and Fabrizio [1] due to the necessity of using of a model describing the behavior of classical viscoelastic materials, thermal media, electromagnetic systems, etc.
 
In [3], Nieto and Losado studied the equation
$$
_{CF}D_{0t}^\alpha f(t)=u(t)
$$
and based on its solution, they have introduced an integral operator corresponding to the differential operator (1) as
$$
_{CF}I_{0t}^\alpha f(t)=(1-\alpha)u(t)+\alpha\int\limits_0^t u(s)ds,\,\,\,t\geq 0.
$$

In their work, they have also considered the following initial value problem
$$
\begin{array}{l}
_{CF}D_{0t}^\alpha f(t)=\lambda f(t)+u(t),\,\,t\geq 0,\\
f(0)=f_0\in \mathbb{R}.\\
\end{array}
$$
This problem has been reduced to a first order ODE and they have proved that the problem has a unique solution for any $\lambda\in \mathbb{R}$.

In this work, we present explicit forms of the solution to the same problem, imposing the required conditions to the given data, by reducing it to the second kind Volterra integral equation.

Another application of this new derivative was considered by Atangana [4], where he studied using nonlinear Fisher's reaction-diffusion equation by using Sumudu transform.

We would like also to note that there is another new fractional order operator without singular kernel, which is an analog of the Riemann-Liouville fractional derivative with singular kernel, was proposed by Yang et al [5] along with its application to the steady heat flow process.

In the next section, we will present our results regarding a unique solvability of certain initial value problem. The main result of this work is presented in section 2, where four different boundary value problems were considered. 
  
\subsection{Solution of initial value problem}

Here we will consider the following initial value problem:

\textbf{IVP.} Find a solution of the equation
\begin{equation}\label{eq2}
_{CF}D_{0t}^\alpha u(t)-\lambda u(t)=f(t),\,\,\,0\leq t\leq T,
\end{equation}
satisfying the initial condition
\begin{equation}\label{eq3}
u(0)=0,
\end{equation}
where $f(t)$ is a given function and $\lambda,\,\alpha\in \mathbb{R}$ such that $0< \alpha\leq 1$.

A unique solvability of this problem is formulated in the following theorem.

\textbf{Theorem 1.} 
\begin{itemize}
\item If $\lambda\neq \frac{1}{1-\alpha}$, $f(t)\in C(0,\infty)$ and $f(0)=0$, then problem (2)-(3) has a unique continuous solution, which is given by
\begin{equation}\label{eq4}
u(t)=\frac{1-\alpha}{1-\lambda(1-\alpha)}f(t)+\frac{\alpha}{\left[1-\lambda(1-\alpha)\right]^2}\int\limits_0^t{f(\xi)e^{\frac{\lambda\alpha}{1-\lambda(1-\alpha)}(t-\xi)}d\xi};
\end{equation}

\item If $\lambda=\frac{1}{1-\alpha}$, $f(t)\in C^1(0,\infty)$ and $f(0)=0,\,f'(0)=0$, then a unique continuous solution to problem (2)-(3) exists and given by
\begin{equation}\label{eq5}
u(t)=-(1-\alpha)f(t)-\frac{(1-\alpha)^2}{\alpha}f'(t).
\end{equation}
\end{itemize}

\begin{proof}

Using the definition of Caputo-Fabrizio operator and integrating by parts, we deduce
\begin{equation}\label{eq6}
\left(\frac{1}{1-\alpha}-\lambda\right)u(t)-\frac{\alpha}{(1-\alpha)^2}\int\limits_0^t{u(s)e^{-\frac{\alpha}{1-\alpha}(t-s)}ds}=f(t).
\end{equation}

First, we consider the case of $\lambda\neq \frac{1}{1-\alpha}$ and hence, (6) yields 
\begin{equation}\label{eq7}
u(t)-\int\limits_0^t{u(s)K(t,s)ds}=\overline{f}(t),
\end{equation}
where
$$
K(t,s)=\frac{\alpha}{(1-\alpha)\left[1-\lambda(1-\alpha)\right]}e^{-\frac{\alpha}{1-\alpha}(t-s)},\,\,\,\,\,
\overline{f}(t)=\frac{1-\alpha}{1-\lambda(1-\alpha)}f(t).
$$
The second kind Volterra integral equation (7) can be solved by the method of successive iterations as follows:
$$
u_1(t)=\overline{f}(t)+\int\limits_0^t{\overline{f}(s)K(t,s)ds}.
$$
$$
\begin{array}{l}
u_2(t)=\overline{f}(t)+\int\limits_0^t{u_1(s)K(t,s)ds}=\overline{f}(t)+\int\limits_0^t{\overline{f}(s)K(t,s)ds}+\\
+\int\limits_0^t{\overline{f}(\xi)d\xi\int\limits_\xi^t{K(t,s)K(s,\xi)ds}}.\\
\end{array}
$$
Setting 
$$
K_2(t,\xi)=\int\limits_\xi^t{K(t,s)K(s,\xi)ds},
$$
$u_2(t)$ can be rewritten as
$$
u_2(t)=\overline{f}(t)+\int\limits_0^t{\overline{f}(\xi)\left[K(t,\xi)+K_2(t,\xi)\right]d\xi}.
$$
One can then prove by mathematical induction that
$$
u_n(t)=\overline{f}(t)+\int\limits_0^t{\overline{f}(\xi)\sum\limits_{i=1}^{n}K_i(t,\xi)d\xi},
$$
where
$$
K_1(t,\xi)=K(t,\xi),\,\,\,K_j(t,\xi)=\int\limits_\xi^t{K(t,s)K_{j-1}(s,\xi)ds},\,\,j=2,3,...
$$

Similarly, we can find general expression for kernels $K_j(t,\xi)$. 
$$
\begin{array}{l}
K_2(t,\xi)=\int\limits_\xi^t{K(t,s)K_1(s,\xi)ds}=\\
=\int\limits_\xi^t{\frac{\alpha}{(1-\alpha)\left[1-\lambda(1-\alpha)\right]}e^{-\frac{\alpha}{1-\alpha}(t-s)}\times \frac{\alpha}{(1-\alpha)\left[1-\lambda(1-\alpha)\right]}e^{-\frac{\alpha}{1-\alpha}(s-\xi)}ds}=\\
=\left(\frac{\alpha}{(1-\alpha)\left[1-\lambda(1-\alpha)\right]}\right)^2\int\limits_\xi^t{e^{-\frac{\alpha}{1-\alpha}(t-\xi)}ds}=\left(\frac{\alpha}{(1-\alpha)\left[1-\lambda(1-\alpha)\right]}\right)^2(t-\xi)e^{-\frac{\alpha}{1-\alpha}(t-\xi)},\\
\end{array}
$$
$$
\begin{array}{l}
K_3(t,\xi)=\int\limits_\xi^t{K(t,s)K_2(s,\xi)ds}=\\
=\int\limits_\xi^t{\frac{\alpha}{(1-\alpha)\left[1-\lambda(1-\alpha)\right]}e^{-\frac{\alpha}{1-\alpha}(t-s)}\times \left(\frac{\alpha}{(1-\alpha)\left[1-\lambda(1-\alpha)\right]}\right)^2(s-\xi)e^{-\frac{\alpha}{1-\alpha}(s-\xi)}ds}=\\
=\left(\frac{\alpha}{(1-\alpha)\left[1-\lambda(1-\alpha)\right]}\right)^3\frac{(t-\xi)^2}{2}e^{-\frac{\alpha}{1-\alpha}(t-\xi)},\\
\end{array}
$$
$\ldots\,\,\,\,\,\ldots$
$$
K_i(t,\xi)=\left(\frac{\alpha}{(1-\alpha)\left[1-\lambda(1-\alpha)\right]}\right)^i\frac{(t-\xi)^{i-1}}{(i-1)!}e^{-\frac{\alpha}{1-\alpha}(t-\xi)},\,\,\,\,i=1,2,...
$$
Hence, resolvent-kernel will have the form
$$
\begin{array}{l}
\displaystyle{R(t,\xi)=\sum\limits_{i=1}^{\infty}K_i(t,\xi)
=\frac{\alpha}{(1-\alpha)\left[1-\lambda(1-\alpha)\right]}e^{-\frac{\alpha}{1-\alpha}(t-\xi)}}\times\\
\displaystyle{\times\sum\limits_{i=1}^{\infty}\frac{\left[\frac{\alpha}{(1-\alpha)\left[1-\lambda(1-\alpha)\right]}(t-\xi)\right]^{i-1}}{(i-1)!},}\\
\end{array}
$$
which can be reduced to
$$
R(t,\xi)=\frac{\alpha}{(1-\alpha)\left[1-\lambda(1-\alpha)\right]}e^{\frac{\lambda\alpha}{1-\lambda(1-\alpha)}(t-\xi)}.
$$
Thus, solution of (7) will be given by
$$
u(t)=\overline{f}(t)+\frac{\alpha}{(1-\alpha)\left[1-\lambda(1-\alpha)\right]}\int\limits_0^t{\overline{f}(\xi)e^{\frac{\lambda\alpha}{1-\lambda(1-\alpha)}(t-\xi)}d\xi},
$$
which on using the representation of $\overline{f}(t)$ and the condition $f(0)=0$ leads to the solution representation (4) as desired.

Now, if $\lambda=\frac{1}{1-\alpha}$, equation (6) reduces to
$$
\int\limits_0^t{u(s)e^{-\frac{\alpha}{1-\alpha}(t-s)}ds}=\frac{(1-\alpha)^2}{\alpha}f(t),
$$
which can be rewritten as the following second kind Volterra integral equation
\begin{equation}\label{eq8}
u(t)-\int\limits_0^t u(s)\overline{K}(t,s)ds=\widehat{f}(t),
\end{equation}
where
$$
\overline{K}(t,s)=\frac{\alpha}{1-\alpha}e^{-\frac{\alpha}{1-\alpha}(t-s)},\,\,\,\,
\widehat{f}(t)=\frac{(1-\alpha)^2}{\alpha}f'(t).
$$
Following the same previous approach, we obtain an expression for the resolvent-kernel in the form $\overline{R}(t,s)=\frac{\alpha}{1-\alpha}$ and hence, a solution of integral equation (8) will be given by
$$
u(t)=\widehat{f}(t)+\frac{\alpha}{1-\alpha}\int\limits_0^t\widehat{f}(s)ds,
$$
which on using the representation of $\widehat{f}$ and conditions $f(0)=0,\,\,f'(0)=0$, leads to the explicit form of the solution as in (5). This ends the proof of Theorem 1.
\end{proof}

\textbf{Remark 1.} The solution for the special case $\lambda=0$ follows from the case $\lambda\neq \frac{1}{1-\alpha}$ and is given by
$$
u(t)=(1-\alpha)f(t)+\alpha\int\limits_0^t f(s)ds.
$$  

\textbf{Remark 2.} If equation (2) is subjected to a non-homogeneous initial condition $u(0)=u_0$, then the condition $f(0)=0$ will be replaced by $f(0)=-\lambda u_0$ and hence the solution will be given by
$$
u(t)=-\frac{(1-\alpha)^2}{\alpha}f'(t)-(1-\alpha)[f(t)-f(0)]+u_0\quad \mbox{for} \quad \lambda=\frac{1}{1-\alpha} 
$$
and 
$$
\begin{array}{l}
\displaystyle{u(t)=\frac{1-\alpha}{1-\lambda(1-\alpha)}f(t)+\frac{\alpha}{\left[1-\lambda(1-\alpha)\right]^2}\int\limits_0^t{f(\xi)e^{\frac{\lambda\alpha}{1-\lambda(1-\alpha)}(t-\xi)}d\xi}+}\\
\displaystyle{+\frac{u_0}{1-\lambda(1-\alpha)}e^{\frac{\lambda\alpha}{1-\lambda(1-\alpha)}t} \quad \quad \mbox{for}\quad \lambda\neq\frac{1}{1-\alpha}.}\\
\end{array}
$$

\section{Main result}

In this section we will consider Caputo-Fabrizio fractional heat equation subjected to four different boundary conditions associated with self-adjoint and non self-adjoint spectral problems.

\subsection{Boundary value problems associated with self-adjoint spectral problems}

Consider a rectangular domain $\Omega=\left\{(x,t):\,0<x<1,\,0<t<T\right\}$. In this domain we investigate the following three problems:

\textbf{Problem 1.} Find a regular solution of the equation
\begin{equation}\label{eq9}
{}_{CF}D_{0t}^\alpha u(x,t)-u_{xx}(x,t)=g(x,t),
\end{equation}
in a domain $\Omega$, which satisfies initial condition
\begin{equation}\label{eq10}
u(x,0)=0,\,\,\,0\leq x\leq 1
\end{equation}
and boundary conditions
\begin{equation}\label{eq11}
u(0,t)=0,\,\,\,\,\,u(1,t)=0, \,\,0\leq t\leq T,
\end{equation}
where $g(x,t)$ is a given function.

\textbf{Problem 2.} Find a regular solution of problem (9)-(10) in the domain $\Omega$, which satisfies boundary conditions
\begin{equation}\label{eq12}
u'(0,t)=0,\,\,\,\,\,u'(1,t)=0, \,\,0< t< T.
\end{equation}

\textbf{Problem 3.} Find a regular solution of problem (9)-(10) in $\Omega$, which satisfies non-local boundary conditions
\begin{equation}\label{eq13}
u(0,t)=u(1,t),\,\,0\leq t\leq T;\,\,u'(0,t)=u'(1,t), \,\,0< t< T.
\end{equation}

We will first consider Problem 1. Using the method of separation variables leads to the following self-adjoint spectral problem
$$
X''(x)+\mu X(x)=0,\,\,\,X(0)=X(1)=0
$$
in the variable $x$. This problem has eigenvalues $\mu_k=-(k\pi)^2,\,k=1,2,...$ and the corresponding eigenfunctions are $X_k(x)=\sin k\pi x$.

Since the system of functions $\{\sin k\pi x\}$ is complete and forms a basis in $L_2$, we look for a solution to Problem 1 of the form
\begin{equation}\label{eq14}
u(x,t)=\sum\limits_{k=1}^\infty u_k(t)\sin k\pi x,
\end{equation} 
where
\begin{equation}\label{eq15}
u_k(t)=\int\limits_0^1{u(x,t)\sin k\pi x dx}.
\end{equation}
Substituting (14) into (9) and (10), we get
\begin{equation}\label{eq16}
{}_{CF}D_{0t}^\alpha u_k(t)+(k\pi)^2 u_k(t)=g_k(t),\,\,\,\,u_k(0)=0
\end{equation} 
where
\begin{equation}\label{eq17}
g_k(t)=\int\limits_0^1{g(x,t)\sin k\pi x dx}.
\end{equation}
According to Theorem 1, solution of problem (16) is given by
\begin{equation}\label{eq18}
u_k(t)=\frac{(1-\alpha)g_k(t)}{1+(k\pi)^2(1-\alpha)}+\frac{\alpha}{\left[1+(k\pi)^2(1-\alpha)\right]^2}\int\limits_0^tg_k(\xi)e^{-\frac{\alpha(k\pi)^2(t-\xi)}{1+(k\pi)^2(1-\alpha)}}d\xi
\end{equation}
with $g_k(t)\in C[0,T],\,g_k(0)=0$, which can be achieved by assuming $g(x,t)\in C[0,T]$ and $g(x,0)=0$.
This imposed conditions will in turn, lead to the convergence of the series solution given by (14). 

Series representation of $u_{xx}(x,t)$ is given by
$$
\begin{array}{l}
u_{xx}(x,t)=-\sum\limits_{k=1}^\infty (k\pi)^2\frac{1-\alpha}{1+(k\pi)^2(1-\alpha)}g_k(t)\sin k\pi x-\\
-\sum\limits_{k=1}^\infty (k\pi)^2\left[\frac{\alpha}{\left[1+(k\pi)^2(1-\alpha)\right]^2}\int\limits_0^tg_k(\xi)e^{-\frac{\alpha(k\pi)^2(t-\xi)}{1+(k\pi)^2(1-\alpha)}}d\xi\right]\sin k\pi x=\\
=\sum\limits_{k=1}^\infty g_k(t)\sin k\pi x +\sum\limits_{k=1}^\infty\frac{\sin k\pi x}{1+(k\pi)^2(1-\alpha)}\int\limits_0^t g_k'(\xi)e^{-\frac{\alpha(k\pi)^2(t-\xi)}{1+(k\pi)^2(1-\alpha)}}d\xi.
\end{array}
$$

The convergence of the second series in the expression of $u_{xx}(x,t)$ is guaranteed by assuming $g_t(x,t)\in L_1[0,T]$.  Hence, in order to prove the convergence of series expansion of $u_{xx}(x,t)$ it remains to show the convergence of the first series in the expression of $u_{xx}(x,t)$. To do so, we consider the following estimate of this series: 
$$
\begin{array}{l}
\sum\limits_{k=1}^\infty |\int\limits_0^1 g(x,t)\sin k\pi x dx|=\sum\limits_{k=1}^\infty |\frac{1}{k\pi}\left(g(1,t)(-1)^k-g(0,t)\right)+\\
+\frac{1}{k\pi}\int\limits_0^1 \frac{\partial g(x,t)}{\partial x}\sin k\pi x dx|,\\
\end{array}
$$
which on assuming $g(0,t)=g(1,t)=0$ and using the inequality $ab\leq 1/2(a^2+b^2)$ becomes
$$
\sum\limits_{k=1}^\infty \left|\int\limits_0^1 g(x,t)\sin k\pi x dx\right|\leq \sum\limits_{k=1}^\infty\frac{1}{2}\left(\frac{1}{(k\pi)^2}+|\overline{g_k}(t)|^2\right),
$$
where
$$
\overline{g_k}(t)=\int\limits_0^1 \frac{\partial g(x,t)}{\partial x}\sin k\pi x dx.
$$

The convergence of the latter series is obtained by assuming $\frac{\partial g(x,t)}{\partial x}\in L_2[0,1]$ and using $\sum\limits_{k=1}^\infty \overline{g_k}(t)\leq \Vert\overline{g_k}(t)\Vert_{L_2}$. Hence, series representation for $u_{xx}(x,t)$ converges. The convergence of series expansion for ${}_CD_{0t}^\alpha u(x,t)$ follows from equation (9).

This result can be formulated in the following theorem.

\textbf{Theorem 2.} If the following conditions
$$
g(x,t)\in C\left(\overline{\Omega}\right),\, g(x,0)=0,\,\,g(0,t)=g(1,t)=0,\,g_t(x,t)\in L_1[0,T],\,\,g_x(x,t)\in L_2[0,1],
$$
hold, then Problem 1 has a unique solution represented by 
\begin{equation}\label{eq19}
u(x,t)=\sum\limits_{k=1}^\infty\left[\frac{(1-\alpha)g_k(t)}{1+(k\pi)^2(1-\alpha)}+\frac{\alpha\int\limits_0^t g_k(\xi)e^{-\frac{\alpha(k\pi)^2(t-\xi)}{1+(k\pi)^2(1-\alpha)}}d\xi}{\left[1+(k\pi)^2(1-\alpha)\right]^2}\right]\sin k\pi x.
\end{equation}

Note that a uniqueness of the solution for Problem 1 will follow from the representation (15), based on (18) and completeness of the system $\{\sin k\pi x\}$.

Since the boundary conditions in Problems 2 and 3 will result in self-adjoint spectral problems, then Problems 2 and 3 could be studied in a similar way. The results can be formulated in the following theorems.

\textbf{Theorem 3.} If the following conditions
$$
g(x,t)\in C(\overline{\Omega}),\,\,g_t(x,t)\in L_1[0,T],\,\,g_x(x,t)\in L_2[0,1]
$$
hold, then a unique solution of Problems 2, represented by
$$
u(x,t)=\sum\limits_{n=0}^\infty\left[\frac{(1-\alpha)g_n(t)}{1+(n\pi)^2(1-\alpha)}+\frac{\alpha\int\limits_0^t g_n(\xi)e^{-\frac{\alpha(n\pi)^2(t-\xi)}{1+(n\pi)^2(1-\alpha)}}d\xi}{\left[1+(n\pi)^2(1-\alpha)\right]^2}\right]\cos n\pi x
$$
exists, where
$$
g_n(t)=\int\limits_0^1 g(x,t)\cos n\pi x dx,\,\,n=0,1,2,...
$$

\textbf{Theorem 4.} If the following conditions
$$
g(x,t)\in C(\overline{\Omega}),\,\,g(0,t)=g(1,t),\,\,g_t(x,t)\in L_1[0,T],\,\,g_x(x,t)\in L_2[0,1]
$$
hold, then a unique solution of Problems 3, represented by
$$
\begin{array}{l}
u(x,t)=\\
\displaystyle{=\sum\limits_{n=0}^\infty\left[\frac{(1-\alpha)g_{n1}(t)}{1+(n\pi)^2(1-\alpha)}+\frac{\alpha\int\limits_0^t g_{n1}(\xi)e^{-\frac{\alpha(n\pi)^2(t-\xi)}{1+(n\pi)^2(1-\alpha)}}d\xi}{\left[1+(n\pi)^2(1-\alpha)\right]^2}\right]\cos 2n\pi x}+\\
+
\displaystyle{\sum\limits_{n=1}^\infty\left[\frac{(1-\alpha)g_{n2}(t)}{1+(n\pi)^2(1-\alpha)}+\frac{\alpha\int\limits_0^t g_{n2}(\xi)e^{-\frac{\alpha(n\pi)^2(t-\xi)}{1+(n\pi)^2(1-\alpha)}}d\xi}{\left[1+(n\pi)^2(1-\alpha)\right]^2}\right]\sin 2n\pi x}\\
\end{array}
$$
exists, where
$$
g_{n1}(t)=\int\limits_0^1 g(x,t)\cos 2n\pi x dx,\,\,n=0,1,2,...,
$$
$$
g_{n2}(t)=\int\limits_0^1 g(x,t)\sin 2n\pi x dx,\,\,n=1,2,...
$$

\subsection{Boundary value problem associated with nonself-adjoint spectral problem}

In this section, we consider the following problem with non-local boundary conditions: 

\textbf{Problem 4. } Find a regular solution to problem (9)-(10) in $\Omega$, which satisfies the non-local boundary conditions
\begin{equation}\label{eq20}
u(0,t)=u(1,t),\,\,0\leq t\leq T,\,\,\,\,u_x(0,t)=0, \,\,0< t< T.
\end{equation}

The associated spectral problem for this case is given by
\begin{equation}\label{eq21}
{X}''\left( x \right)+\mu X\left( x \right)=0,\,\,X\left( 0 \right)=X\left( 1 \right),\,\,{X}'\left( 0 \right)=0,
\end{equation}
which is not self-adjoint. The eigenvalues of (21) are ${{\mu }_{k}}=\lambda _{k}^{2},\,\,{{\lambda }_{k}}=2k\pi \,(k=0,1,2,...)$ and the corresponding eigenfunctions are $1,\,\cos {{\lambda }_{k}}x$, supplemented by the associate function $x\sin {{\lambda }_{k}}x$, which form a complete system of root functions denoted by
\begin{equation}\label{eq22}
{{X}_{k}}\left( x \right)=\left\{ 1,\,\cos {{\lambda }_{k}}x,\,x\sin {{\lambda }_{k}}x \right\},\,\,k=1,2,....
\end{equation}

Since, problem (21) is not self-adjoint and hence $X_k$ does not form a basis, we need to find root functions of the corresponding adjoint problem:
$$
{Y}''\left( x \right)+\mu Y\left( x \right)=0,\,\,{Y}'\left( 0 \right)={Y}'\left( 1 \right),\,\,Y\left( 1 \right)=0.
$$
This problem has the following system of root functions:
\begin{equation}\label{eq23}
{{Y}_{k}}\left( x \right)=\left\{ 2(1-x),4(1-x)\,\cos {{\lambda }_{k}}x,\,4\sin {{\lambda }_{k}}x \right\}\,\,k=1,2,....
\end{equation}

Now systems (22) and (23) form bi-orthogonal system, which satisfies the necessary and sufficient condition for the basis property in the space ${{L}_{2}}[0,1]$ (see [6]).

Thus, we seek a solution of Problem 4 in the form
\begin{equation}\label{eq24}
u(x,t)={{u}_{0}}(t)+\sum\limits_{k=1}^{\infty }{}{{u}_{1k}}(t)\cos 2k\pi x+\sum\limits_{k=1}^{\infty }{}{{u}_{2k}}(t)\, x\sin 2k\pi x,\,\,0\leq t\leq T.
\end{equation}
The given function $g(x,t)$ can be also represented in the following series expansion
\begin{equation}\label{eq25}
g(x,t)={{g}_{0}(t)}+\sum\limits_{k=1}^{\infty }{}{{g}_{1k}(t)}\cos 2k\pi x+\sum\limits_{k=1}^{\infty }{}{{g}_{2k}(t)}\, x\sin 2k\pi x,\,\,\,0\leq t\leq T,
\end{equation}
where the coefficients of the two series above are defined as follows
\begin{equation}\label{eq26}
\begin{array}{l}
{{u}_{0}}(t)=2\int\limits_{0}^{1}{}u(x,t)(1-x)\,dx,\\
{{u}_{1k}}(t)=4\int\limits_{0}^{1}{}u(x,t)(1-x)\cos 2k\pi x\,dx,\\
{{u}_{2k}}(t)=4\int\limits_{0}^{1}{}u(x,t)\sin 2k\pi x\,dx,\\
{{g}_{0}}(t)=2\int\limits_{0}^{1}{}g(x,t)(1-x)\,dx,\\
{{g}_{1k}}(t)=4\int\limits_{0}^{1}{}g(x,t)(1-x)\cos 2k\pi x\,dx,\\
{{g}_{2k}}(t)=4\int\limits_{0}^{1}{}g(x,t)\sin 2k\pi x\,dx. \\
\end{array}
\end{equation}

The unknown coefficients $u_0(t),\,u_{1k}(t),\,u_{2k}(t)$ can be determined by substituting (23)-(24) into (9) and (10), as solution of the following fractional order initial value problems:
\begin{equation}\label{eq27}
_{CF}D_{0t}^{\alpha }{{u}_{0}}\left( t \right)={{g}_{0}}(t),\,\,u_0(0)=0,
\end{equation}
\begin{equation}\label{eq28}
_{CF}D_{0t}^{\alpha }{{u}_{1k}}\left( t \right)+{{\left( 2k\pi  \right)}^{2}}{{u}_{1k}}\left( t \right)={{g}_{1k}}(t)+4k\pi {{u}_{2k}}\left( t \right),\,\,u_{1k}(0)=0,
\end{equation}
\begin{equation}\label{eq29}
_{CF}D_{0t}^{\alpha }{{u}_{2k}}\left( t \right)+{{\left( 2k\pi  \right)}^{2}}{{u}_{2k}}\left( t \right)={{g}_{2k}}(t),\,\,u_{2k}(0)=0.
\end{equation}

Based on Theorem 1, solutions of (27)  and (29) are given, respectively, by
\begin{equation}\label{eq30}
u_0(t)=(1-\alpha)g_0(t)+\alpha\int\limits_0^t g_0(z)dz,
\end{equation}
\begin{equation}\label{eq31}
\begin{array}{l}
u_{2k}(t)=\frac{1-\alpha}{1+(2k\pi)^2(1-\alpha)}g_{2k}(t)+\\
+\frac{\alpha}{\left[1+(2k\pi)^2(1-\alpha)\right]^2}\int\limits_0^t g_{2k}(\xi)e^{-\frac{\alpha(2k\pi)^2(t-\xi)}{1+(2k\pi)^2(1-\alpha)}}d\xi.\\
\end{array}
\end{equation}

Similarly, on using the expression (31), solution of (28) can be written as
\begin{equation}\label{eq32}
\begin{array}{l}
u_{1k}(t)=\frac{1-\alpha}{1+(2k\pi)^2(1-\alpha)}\left[ g_{1k}(t)+\frac{4k\pi(1-\alpha)}{1+(2k\pi)^2(1-\alpha)}g_{2k}(t)\right]+\\
+\frac{\alpha}{\left[1+(2k\pi)^2(1-\alpha)\right]^2}\int\limits_0^t e^{-\frac{\alpha(2k\pi)^2(t-z)}{1+(2k\pi)^2(1-\alpha)}}\times\\
\times\left[g_{1k}(z)+\frac{8k\pi(1-\alpha)}{1+(2k\pi)^2(1-\alpha)}g_{2k}(z)+\frac{4k\pi\alpha}{\left[1+(2k\pi)^2(1-\alpha)\right]^2}g_{2k}(z)(t-z)\right]dz.\\
\end{array}
\end{equation}

This completes the existence of formal solution to Problem 4 as given by (24).

It is now left to check the convergence of the series appeared in $u(x,t)$, $u_{xx}(x,t)$ and ${}_CD_{0t}^\alpha u(x,t)$.

Here we present the convergence of series representation of $u_{xx}(x,t)$ and the rest can be treated similarly.

Using the representation of $u(x,t)$ as given in (24) together with (30)-(32), we obtain the following expression for $u_{xx}(x,t)$:  
\begin{equation}\label{eq33}
\begin{array}{l}
u_{xx}(x,t)=-\sum\limits_{k=1}^\infty \frac{(2k\pi)^2(1-\alpha)\cos2k\pi x}{1+(2k\pi)^2(1-\alpha)}g_{1k}(t)-\\
-\sum\limits_{k=1}^\infty \frac{2(2k\pi)^2(1-\alpha)\cos2k\pi x}{\left[1+(2k\pi)^2(1-\alpha)\right]^2} g_{2k}(t)+\\
+\sum\limits_{k=1}^\infty \frac{(2k\pi)^2(1-\alpha)\cos2k\pi x}{\left[1+(2k\pi)^2(1-\alpha)\right]^2}\int\limits_0^t g_{1k}(z)e^{-\frac{\alpha(2k\pi)^2(t-z)}{1+(2k\pi)^2(1-\alpha)}}dz+\\
+\sum\limits_{k=1}^\infty \frac{4(2k\pi)^2\alpha(1-\alpha)\cos2k\pi x}{\left[1+(2k\pi)^2(1-\alpha)\right]^3}\int\limits_0^t g_{2k}(z)e^{-\frac{\alpha(2k\pi)^2(t-z)}{1+(2k\pi)^2(1-\alpha)}}dz+\\
+\sum\limits_{k=1}^\infty \frac{2(2k\pi)^3\alpha^2\cos2k\pi x}{\left[1+(2k\pi)^2(1-\alpha)\right]^4}\int\limits_0^t g_{2k}(z)(1-z)e^{-\frac{\alpha(2k\pi)^2(t-z)}{1+(2k\pi)^2(1-\alpha)}}dz+\\
+\sum\limits_{k=1}^\infty \frac{(2k\pi)^2(1-\alpha)x\sin2k\pi x}{1+(2k\pi)^2(1-\alpha)}g_{2k}(t)+\\
+\sum\limits_{k=1}^\infty \frac{(2k\pi)^2\alpha x\sin2k\pi x}{\left[1+(2k\pi)^2(1-\alpha)\right]^2}\int\limits_0^t g_{2k}(z)e^{-\frac{\alpha(2k\pi)^2(t-z)}{1+(2k\pi)^2(1-\alpha)}}dz.\\
\end{array}
\end{equation} 

On integration by parts, using inequalities $ab\leq 1/2(a^2+b^2)$,  $\sum\limits_{k=1}^\infty f_k(t)\leq \Vert f_k(t)\Vert_{L_2}$, and imposing the following conditions
$$
g(x,t)\in C(\overline{\Omega}),\, g(x,0)=0,\,g(0,t)=g(1,t),\,g_t(x,t)\in L_1[0,T],\,g_x(x,t)\in L_2[0,1],
$$
we get the estimate for $u_{xx}(x,t)$ as:
$$
|u_{xx}(x,t)|\leq \sum\limits_{k=1}^\infty \frac{C_3}{(k\pi)^2}+\left(\Vert{\bar{g_{1k}}(t)}\Vert_{L_2}+\Vert{g_{2k}(t)}\Vert_{L_2}+\Vert{\bar{g_{2k}}(t)}\Vert_{L_2}\right),
$$
where $C_3$ is a positive constant and
$$
\begin{array}{l}
\bar{g_{1k}}(t)=\int\limits_0^1 \left[g_x(x,t)(1-x)-g(x,t)\right]\sin2k\pi x dx,\\
\bar{g_{2k}}(t)=\int\limits_0^1 g_x(x,t)\cos 2k\pi x dx,\\
\end{array}
$$

This estimate will ensure the convergence of series representation of $u_{xx}(x,t)$.

Uniqueness follows from the representation of solution and the completeness of the used bi-orthogonal system.
 
Now, we can formulate our result as the following

\textbf{Theorem 5.} If 
$$
g(x,t)\in C(\overline{\Omega}),\, g(x,0)=0,\,g(0,t)=g(1,t),\,g_t(x,t)\in L_1[0,T],\,g_x(x,t)\in L_2[0,1],
$$
then Problem 4 has a unique solution, represented by 
$$
\begin{array}{l}
\displaystyle{u(x,t)=(1-\alpha)g_0(t)+\alpha\int\limits_0^t g_0(z)dz+}\\
+\displaystyle{\sum\limits_{k=1}^{\infty }\left(\frac{1-\alpha}{1+(2k\pi)^2(1-\alpha)}\left[ g_{1k}(t)+\frac{4k\pi(1-\alpha)}{1+(2k\pi)^2(1-\alpha)}g_{2k}(t)\right]+\right.}\\
\displaystyle{+\frac{\alpha}{\left[1+(2k\pi)^2(1-\alpha)\right]^2}\int\limits_0^t e^{-\frac{\alpha(2k\pi)^2(t-z)}{1+(2k\pi)^2(1-\alpha)}}\left[ g_{1k}(z)+\frac{8k\pi(1-\alpha)}{1+(2k\pi)^2(1-\alpha)}g_{2k}(z)+\right.}\\
\left.\displaystyle{\left.+\frac{4k\pi\alpha}{\left[1+(2k\pi)^2(1-\alpha)\right]^2}g_{2k}(z)(t-z)\right]dz}\right)\cos 2k\pi x+\\
\displaystyle{+\sum\limits_{k=1}^{\infty }\left(\frac{1-\alpha}{1+(2k\pi)^2(1-\alpha)}g_{2k}(t)+\right.}\\
\left.\displaystyle{+\frac{\alpha}{\left[1+(2k\pi)^2(1-\alpha)\right]^2}\int\limits_0^t g_{2k}(\xi)e^{-\frac{\alpha(2k\pi)^2(t-\xi)}{1+(2k\pi)^2(1-\alpha)}}d\xi}\right) x\sin 2k\pi x.\\
\end{array}
$$
\section{Acknowledgement} 
Authors acknowledge  financial support from The Research Council (TRC), Oman. This work is funded by TRC under the research agreement no. ORG/SQU/CBS/13/030

\end{document}